%
%
\input oldamstex
\documentstyle {oldamsppt}
\magnification=1200
\tolerance=3000
\openup 6 pt
\nologo
\def\SBIMSMark#1#2#3{
 \font\SBF=cmss10 at 10 true pt
 \font\SBI=cmssi10 at 10 true pt
 \setbox0=\hbox{\SBF Stony Brook IMS Preprint \##1}
 \setbox2=\hbox to \wd0{\hfil \SBI #2}
 \setbox4=\hbox to \wd0{\hfil \SBI #3}
 \setbox6=\hbox to \wd0{\hss
             \vbox{\hsize=\wd0 \parskip=0pt \baselineskip=10 true pt
                   \copy0 \break%
                   \copy2 \break%
                   \copy4 \break}}
 \dimen0=\ht6   \advance\dimen0 by \vsize \advance\dimen0 by 8 true pt
                \advance\dimen0 by -\pagetotal
 \dimen2=\hsize \advance\dimen2 by .25 true in
%
%
  \openin2=publishd.tex
  \ifeof2\setbox0=\hbox to 0pt{}
  \else 
     \setbox0=\hbox to 3.1 true in{
                \vbox to \ht6{\hsize=3 true in \parskip=0pt  \noindent  
                {\SBI Published in modified form:}\hfil\break
                \input publishd.tex 
                \vfill}}
  \fi
  \closein2
  \ht0=0pt \dp0=0pt
 \ht6=0pt \dp6=0pt
 \setbox8=\vbox to \dimen0{\vfill \hbox to \dimen2{\copy0 \hss \copy6}}
 \ht8=0pt \dp8=0pt \wd8=0pt
 \copy8
 \message{*** Stony Brook IMS Preprint #1, #2 ***}
}

\SBIMSMark{1996/6}{June 1996}{}
\bigskip
\bigskip
\bigskip
\bigskip

\topmatter
\title
Invariant Measures for Typical Quadratic Maps
\endtitle
\author
 Marco Martens \footnote {Institute for Mathematical Sciences, SUNY at 
Stony Brook, Stony Brook, NY 11794-3651.} 
Tomasz Nowicki \footnote{University of Warsaw, Poland}
\endauthor
\date{June 28 , 1996}
\endtopmatter

\bigskip
\bigskip
\bigskip
\bigskip
\bigskip
\bigskip
\centerline{\bf Abstract.} 

\bigskip

\flushpar
A sufficient geometrical condition for the existence
of absolutely continuous invariant probability measures for $S-$unimodal maps
will be discussed. The Lebesgue typical existence of such 
measures in the quadratic family will be a consequence.

\vskip .5 in

\bigskip
\centerline{\bf 1. Introduction}
\bigskip

\flushpar
A general belief, or hope, in the theory of dynamical systems is that 
typical dynamical systems have well-understood behavior. This belief
has two forms, depending on the meaning of the word ``typical''. It
could refer to the topological generic situation or to the Lebesgue
typical situation in parameter space. In this work {\it typical} will
refer to Lebesgue typical and the behavior of a Lebesgue typical quadratic 
map on the interval will be discussed.

\flushpar
The quadratic family is formed by the maps
$q_t:[-1,1]\to [-1,1]$ with $t\in [0,1]$ and
$$
q_t(x)=-2tx^\alpha+2t-1,
$$
with the critical exponent $\alpha=2$. The maps in this family can be 
classified as follows. The maps in
$$
\Cal{P}=\{ t\in [0,1]| \text{  } q_t \text{ has a periodic attractor}     \}
$$
have a unique periodic orbit whose basin of attraction is an open and dense 
set. Moreover this basin has full Lebesgue measure. In particular the
invariant measure on the periodic attractor is the SBR-measure for the map.
Recall that a measure $\mu$ on $[-1,1]$ is called an 
S(inai)-B(owen)-R(uelle)-measure for $q_t$ if
$$
\lim_{n\to \infty} \frac1n\sum_{k=0}^{n-1} \delta_{q_t^k(x)}=\mu,
$$
for typical $x\in [-1,1]$. 

\flushpar
The maps in
$$
\Cal{R}=\{ t\in [0,1]|  \text{  }q_t \text{ is infinitely renormalizable} \}
$$
have a unique invariant minimal Cantor set which attracts both generic
and typical orbits. This Cantor set is uniquely ergodic and has zero
Lebesgue measure, [BL2],[G], [M]. 
The unique invariant measure on this Cantor set is the 
SBR-measure for the system.
The maps in
$$
\Cal{I}=[0,1]\setminus \{ \Cal{P}\cup\Cal{R} \}
$$
have a periodic interval whose orbit is the limit set of generic orbits. 
The orbit of this periodic interval absorbs also the orbit of
typical points. These maps are ergodic with respect to the Lebesgue measure, 
[BL1], [GJ], [K], [M]. In the quadratic family, $\alpha=2$, the limit set of
typical points is actually also the orbit of this periodic interval,
[L1]. However, in families with $\alpha$ big enough there are maps in 
$\Cal{I}$ whose typical limit set is not this periodic interval, [BKNS].

\bigskip

\flushpar
Before discussing the behavior of typical quadratic maps 
let us include the behavior of generic quadratic maps.

\proclaim{Theorem ([GS], [L3])} Hyperbolicity is dense in the quadratic
family, $\overline{\Cal{P}}=[0,1]$.
\endproclaim

\flushpar
We will continue to specify the behavior of a typical map in $\Cal{I}$.
 The dynamics of maps in
$$
\Cal{M}=\{t\in \Cal{I}|  \text{  }
q_t \text{ has an absolutely continuous invariant
probability measure}\}
$$
is well-understood. The measure is unique and  its support is the orbit of 
the above periodic interval. Moreover it has positive Lyaponov exponent, 
[K], [Ld]. Starting in [NU], where it was shown that $q_1\in\Cal{M}$, more
and more maps $q_t$ were shown to have such a measure ([B], [R], [Mi]). Finally
it was shown in [Ja] that $\Cal{M}$ has positive measure. 

\proclaim{Theorem A (joint with Lyubich)} A typical quadratic map has a 
unique SBR-measure. More specifically
\parindent=20pt
\item{1)} For $t\in \Cal{P}$ the support of the SBR-measure is the periodic 
          attractor,
\item{2)} For $t\in \Cal{R}$ the SBR-measure is supported on a Cantor set.
\item{3)} For $t\in \Cal{M}$ the SBR-measure is an absolutely continuous 
measure supported on the orbit of a periodic interval.
\item{4)} the set $\Cal{P}\cup \Cal{R}\cup \Cal{M}\subset [0,1]$ has full
Lebesgue measure.
\endproclaim

\flushpar
Johnson constructed 
 unimodal maps in $\Cal{I}$ (with arbitrary critical 
exponent) which do not have an
absolutely continuous invariant measure, [Jo]. More careful 
combinatorial Johnson-Examples were made without SBR-measure ([HK]). 
The same work 
shows the existence of maps in $\Cal{I}\setminus \Cal{M}$ which have
an SBR-measure but this measure is not absolutely continuous. The complications
which occur in $\Cal{I}\setminus \Cal{M}$ are thoroughly studied in [Br].

\bigskip

\flushpar
In this work we will formulate a geometrical condition on maps in 
$\Cal{I}$ sufficient for the existence of absolutely continuous invariant 
probability measures. 
The geometric condition is formulated in terms of a decreasing sequence of 
{\it central intervals} $U_n=(-u_n,u_n)$, $n\ge 1$, which are defined for
all unimodal maps with recurrent critical orbit. The domain $D_n\subset U_n$
of the first return map $R_n:D_n\to U_n$, $n\ge 1$ is a countable collection 
of intervals. The central component of $D_n$ is $U_{n+1}$. The first return 
map $R_n$ is said to have a {\it central return} when
$$
R_n(0)\in U_{n+1}.
$$ 
The sufficient geometrical condition for the existence of absolutely
continuous measures is stated in terms of {\it scaling factors}
$$
\sigma_n=\frac{u_{n+1}}{u_n}, n\ge 1.
$$
These scaling factors describe the small scale geometrical properties of the
system but they are also strongly related to distortion questions.
The main consequence of the distortion Theory developed in [M] are the 
a priori bounds
on the distortion of each $R_n$. The renormalization Theory developed in [L1] 
and [LM] achieved much stronger results: if a quadratic unimodal map has only 
finitely many central returns then the scaling factors tend exponentially to 
zero.

\flushpar
The scaling factors are related to small scale geometry, distortion but also 
expansion. The technical step in this work is to show that small 
scaling factors imply strong expansion along the critical orbit. In [NS] 
it was shown that enough expansion along the critical orbit causes the 
existence of an absolutely continuous invariant probability measure.

\proclaim{Theorem B} If an $S-$unimodal map with critical exponent $\alpha>1$
has summable scaling factors, that is
$$
\sum_{n\ge 1} \sigma_n^{1/\alpha} < \infty
$$
then it has an absolutely continuous invariant probability measure.
\endproclaim

\proclaim{Corollary} If a quadratic map has 
only finitely many central returns then it has an absolutely continuous
invariant probability measure.
\endproclaim

\flushpar
The (Johnson-)Examples in [Jo] 
have infinitely many (cascades of) central returns. 
The corollary states that the only quadratic 
unimodal maps in $\Cal{I}$ which do
not have an absolutely continuous invariant  measure are Johnson-Examples. 
The families $\{q_t\}$ with $\alpha$
big enough have maps in $\Cal{I}$ which do not have an absolutely
continuous invariant probability measure and which are also not 
Johnson-Examples, [BKNS].

\flushpar
In [L4], Lyubich studies the parameter space of the (holomorphic)
quadratic family. A new proof showing that $\Cal{I}$ has positive Lebesgue 
measure is given  (compare with the Jacobson-Theorem [Ja]). 
Moreover it is shown that for almost every parameter in $\Cal{I}$
the corresponding quadratic map has only finitely many central returns. 
This, together with Theorem B, implies Theorem A.

\proclaim{Conjecture} A typical map in the family $\{q_t\}$, with critical 
exponent $\alpha>1$, has a unique SBR-measure. More specifically
\parindent=20pt
\item{1)} For $t\in \Cal{P}$ the support of the SBR-measure is the periodic 
          attractor,
\item{2)} For $t\in \Cal{R}$ the SBR-measure is supported on a Cantor set.
\item{3)} For $t\in \Cal{M}$ the SBR-measure is an absolutely continuous 
measure supported on the orbit of a periodic interval.
\item{4)} the set $\Cal{P}\cup \Cal{R}\cup \Cal{M}\subset [0,1]$ has full
Lebesgue measure.
\endproclaim

\flushpar
An appendix is added to collect the standard notions 
and Lemmas in interval dynamics.

\bigskip
\centerline{\bf 2. Central Intervals}
\bigskip

\flushpar
Throughout the following sections
we will fix an $S-$unimodal map $f:[-1,1]\to [-1,1]$ with critical exponent 
$\alpha>1$ and without periodic attractors. Furthermore assume that the 
critical orbit is recurrent. 

\flushpar
The set of nice points is
$$
\Cal{N}=\{x\in [-1,1]| \forall i\ge 0 \text{  } f^i(x)\notin (-|x|,|x|)\}
$$ 
This set is closed and not empty. For example the fixed point of $f$ in 
$(0,1)$ is in $\Cal{N}$. 

\flushpar
For $x\in \Cal{N}$ let $D_x\subset U_x=(-|x|,|x|)$ 
be the set of points whose orbit returns to $U_x$. The first return map to 
$U_x$ is denoted by 
$$
R_x:D_x\to U_x.
$$

\flushpar
The next Lemma is a straightforward consequence of
the fact that 
the boundary of each $U_x$ is formed by  nice points.

\proclaim{Lemma 2.1([M])} For every $x\in \Cal{N}$ there exists a collection of
pairwise dispoint intervals $\Cal{U}_x$ with
\parindent=15pt
\item{1)} $I\subset U_x$ for all $I\in \Cal{U}_x$,
\item{2)} $\bigcup_{I\in \Cal{U}_x} I =D_x$,
\item{3)} if $I\in \Cal{U}_x$ and $0\notin I$ then
$R_x:I\to U_x$  
is monotone and onto,
\item{4)} if $I\in \Cal{U}_x$ and $0\in I$ then 
$
R_x:I\to U_x
$
is $2$ to $1$ onto the image. Moreover $R_x(\partial I)= \{x\} 
\text{ or } \{-x\}$.
\endproclaim

\flushpar
Define the function $\psi:\Cal{N}\to\Cal{N}$ by
$$
\psi(x)=\partial V_x\cap(0,1),
$$
where $V_x\in \Cal{U}_x$ is the {\it central} interval: $0\in V_x$. Say
$R_x|_{V_x}=f^{q_x}$ and observe that
$$
\{f(\psi(x), f^2(\psi(x)), \dots, f^{q_x}(\psi(x))=x\}\cap U_x=\emptyset
$$
which follows from the fact that $R_x:D_x\to U_x$ is the first return map.
In particular $\psi(x)\in \Cal{N}$ and we can consider the 
first return map to $V_x$. It will also satisfy Lemma 2.1.

\flushpar
Choose $u_1\in \Cal{N}$
and consider the sequence $u_{n}=\psi(u_{n-1})$ with $n\ge 1$. Use the 
simplified notation $\Cal{U}_n$ for $\Cal{U}_{u_n}$ and denote the first
return maps by
$$
R_n:D_n\to U_n
$$
instead of $R_{u_n}:D_{u_n}\to U_{u_n}$. All these first return maps 
have the  properties stated in Lemma 2.1. Observe that $|U_n|=2u_n$. 

\bigskip

\flushpar
Let $\sigma_n=\frac{u_{n+1}}{u_n}$, $n\ge 1$. We call $\sigma_n$ the 
{\it scaling factor} of level $n$. We will assume that
$$
\sigma_n\to 0.
$$
However, the main Proposition 3.1, can also be proved by using only an 
a priori bound on the scaling factors. The assumption $\sigma_n\to 0$ will
make the exposition less cumbersome.

\proclaim{Lemma 2.2} 
If $I\in \Cal{U}_n$ and $R_n|_I=f^t$ then there exists an interval 
$J\supset f(I)$ such that 
$$
f^{t-1}:J\to U_{n-1}
$$ 
is monotone onto. In particular all the maps $f^{t-1}: f(I) \to U_n$, 
$I\in \Cal{U}_n$ have uniformly bounded distortion. Moreover these maps will
be essentially linear when $n\to \infty$. 
\endproclaim

\demo{Proof} 
Let $I\in \Cal{U}_n$ with $R_n|I=f^t$ and let 
$J\supset f(I)$ be the maximal interval
on which $f^{t-1}$ is monotone. The maximality implies the existence of $i<t-1$
such that $0\in \partial f^i(J)$. Observe that $f^i(f(I))\cap U_n=\emptyset$,
the first return happens after $t-1>i$ steps. 
So $u_n \text{ (or } -u_n\text{)}\in
f^i(J)$. We observed before that the orbit of $f(u_n)$ never enters $U_{n-1}$,
$f^{t-1}(J)\supset U_{n-1}$.
\hfill\hfill\qed $\,\,$ (Lemma 2.2)
\enddemo

\proclaim{Lemma 2.3} For $\epsilon>0$ there exists $n_0\ge 1$ such that 
the hyperbolic length of any $I\in\Cal{U}_n$ is small,
$$
\text{hyp}(I,U_n)\le \epsilon \text{ and also } 
\frac{|f(I)|}{|f(U_n)|}\le \epsilon,
$$
whenever $n\ge n_0$.
\endproclaim

\demo{Proof} Let $I\in \Cal{U}_n$, say $R_n|I=f^t|I$. 
The previous Lemma states the existence of an interval $J\supset f(I)$ 
such that $f^{t-1}:(J,f(I))\to (U_{n-1},U_n)$ is monotone onto. For $n$ large
we see that $U_n$ is a very small middle interval in $U_{n-1}$, it has a very
small hyperbolic length. Because the map $f$ has negative Schwarzian
derivative we get that $f(I)\subset J$ has a very small hyperbolic length.
Observe that $f^{-1}(J)\subset U_n$. Otherwise the orbit of $u_n$ would pass 
through $U_{m-1}$. However this is impossible: we saw before that the orbit 
of $\psi(x)=u_n$ does not cross $U_x=U_{n-1}$.
The Lemma will be proved by pulling back the pair $(J,f(I))$ one step more.
\hfill\hfill\qed $\,\,$ (Lemma 2.3)
\enddemo

\bigskip
\centerline{\bf 3. Derivatives along Recurrent  Orbits}
\bigskip

\flushpar
Let $\rho_n=\min\{\frac{1}{\sigma_{n-1}},\frac{1}{\sigma_n}\}$.
In this section we will prove 

\proclaim{Proposition 3.1} There exist  $n_0\ge 1$, $\theta<1$ and $C>0$ with 
the following property. If $n\ge n_0$, 
$x\in U_{n+1}$ and $R_n^i(x)\notin U_{n+1}$ for $i\le s$
then
$$
|Df^T(f(x))|\ge C\cdot \rho_n \cdot \theta^{-(s-1)},
$$
where $R^s_n(x)=f^T(x)$.
\endproclaim

\flushpar
In [VT] a similar estimate in the case $s=1$ was obtained for circle 
homeomorphisms with irrational rotation number of bounded type.
The proof of Proposition 3.1 will use the following Lemmas 
and notation. Fix
$x\in U_{n+1}$ according to the Proposition, say $R_n^i(x)=f^{t_i}(x)$ with
$i\le s$.

\proclaim{Lemma 3.2} For each $i\le s$ there exists an interval $S_i\ni f(x)$
such that
$$
f^{t_i-1}:S_i\to U_n
$$ 
is monotone and onto.
\endproclaim

\demo{Proof} Lemma 2.2 applied to $U_{n+1}\in \Cal{U}_n$ states the existence 
of $S_1$. Assume that $S_i\ni x$
exists. Then $f^{t_i-1}: S_i\to U_n$ monotone and onto. Moreover 
$f^{t_i-1}(f(x))\in I_{i+1}\in \Cal{U}_n$. Because 
$f^{t_i-1}(f(x))\notin U_{n+1}$ we have that $I_{i+1}\ne U_{n+1}$ and 
$R_n:I_{i+1}\to U_n$ is monotone and onto. Now let 
$S_{i+1}=f^{-(t_i-1)}(I_{i+1})\cap S_i$.
\hfill\hfill\qed $\,\,$ (Lemma 3.2)
\enddemo

\flushpar
Observe that $f^{t_{i-1}-1}(S_i)=I_i\in \Cal{U}_n$.

\proclaim{Lemma 3.3} There exist $n_0\ge 1$ and $K<\infty$ with the 
following property. If the distortion of
$$
f^{t_i-1}:S_i\to U_n \text{ with } n\ge n_0
$$ 
is bigger than $K$ then
$$
I_i\subset (-\frac34 \cdot u_n, \frac34 \cdot u_n).
$$
\endproclaim

\demo{Proof} Lemma 2.2 states that $f^{t_1}:S_1\to U_n$ has a monotone 
extension up to $U_{n-1}$, the map is essentially linear for big enough $n$. 
Hence $i\ge 2$.
Consider the following decomposition
$$
f^{t_i-1}|S_i=f^{t_{i}-t_{i-1}-1}|f(I_i)\circ 
              f|I_i\circ
              f^{t_{i-1}-1}|S_i.
$$ 
The factor $f^{t_{i-1}-1}|S_i$ has a monotone extension up to 
$U_n$. In particular, for big enough $n$, it
is essentially linear. This is because the image $I_i$
has a small hyperbolic length within $U_n$ (Lemma 2.3). The factor
$f^{t_{i}-t_{i-1}-1}|f(I_i)$
has a monotone extension up to $U_{n-1}$ (Lemma 2.2), it is also essentially
linear. The distortion of $f^{t_i-1}|S_i$ can only be caused by the factor
$f|I_i$, $I_i$ has to be very close to $0$. There is some $n_0$
such that $I_i\subset (-\frac34 \cdot u_n, \frac34\cdot u_n)$, 
whenever $n\ge n_0$.
Here we used Lemma 2.3 which states that $I_i$ has also very small hyperbolic 
length within $U_n$.
\hfill\hfill\qed $\,\,$ (Lemma 3.3)
\enddemo

\proclaim{Lemma 3.4} For any $\theta<1$ there exist $n_0\ge 1$ and 
$C<\infty$ such that 
$$
|S_i|\le C\cdot u_{n+1}^\alpha \cdot \theta^{i-1},
$$
whenever $n\ge n_0$ and $i\ge 2$.
\endproclaim

\demo{Proof} Observe that $f(0)\in S_1\supset S_2\supset\dots\supset S_i$
and $S_i\subset f(U_{n+1})\subset S_1$, $i\ge 2$.
Let $L_j\subset S_1$ be the connected component of $S_1-S_j$ with
$L_j\subset f(U_{n+1})$, $2\le j\le i$.  For $n$ big enough we get 
from the proof of Lemma 3.2 and from Lemma 2.3 that the hyperbolic length of 
$S_j$ within $S_{j-1}$ is very small, $2\le j\le i$. It is easily seen that
this implies 
$$
\aligned
|S_i|&\le C\cdot \theta^{i-1}\cdot |L_i|\\
   &\le C\cdot \theta^{i-1}\cdot |f(U_{n+1})|\\
   &\le C\cdot \theta^{i-1}\cdot u_{n+1}^\alpha.
\endaligned
$$
\hfill\hfill\qed $\,\,$ (Lemma 3.4)
\enddemo

\proclaim{Lemma 3.5} There exist $n_0\ge 1$ and  $C>0$ such that the following 
holds for $n\ge n_0$. If $|R_n(U_{n+1})|=|f^{t_1}(U_{n+1})|\ge 
\frac{1}{10} u_n$ then
$$
|Df^{t_1-1}_{|S_1}|\ge C\cdot \frac{u_n}{u_{n+1}^\alpha}.
$$
If $|R_n(U_{n+1})|<\frac{1}{10} u_n$ then
$$
|Df^{t_1-1}_{|S_1}|\ge C\cdot \frac{u_{n-1}}{u_n^\alpha}.
$$
\endproclaim

\demo{Proof} Consider the map $f^{t_1-1}: S_1\to U_n$. From Lemma 2.2 we know 
that this map has a monotone extension up to $U_{n-1}$. The map is essentially
linear because $u_{n-1}>>u_n$ whenever $n$ is big enough. 
The derivative $|Df^{t_1-1}_{|S_1}|$ is
essentially constant and can be estimated by
$$
|Df^{t_1-1}_{|S_1}|\ge C\frac{|R_n(U_{n+1})|}{|f(U_{n+1})|}
                        \ge C\frac{u_n}{u_{n+1}^\alpha}.
$$
Here we used that $|R_n(U_{n+1})|\ge \frac{1}{10} u_n$.
Now consider the other case: 
$|R_n(U_{n+1})|< \frac{1}{10} u_n$. Let $K\supset f(U_{n+1})$
be the interval which is mapped monotonically onto $U_{n-1}$:
$f^{t_1-1}:K\to U_{n-1}$. Observe that $f^{-1}(K)\subset U_n$. This follows 
from the fact that the orbit of $f(u_n)$ never hits $U_{n-1}$, which was 
observed in section 2. Let $K'\subset K$
be such that $f^{t_1-1}: K'\to (-\frac34 \cdot u_{n-1}, 
                                 \frac34\cdot u_{n-1})$ is monotone and onto.
This map 
has uniform bounded distortion because it has a monotone extension up to
$U_{n-1}$. Let $K''=f^{-1}(K')\subset U_n$. The derivative  
$|Df^{t_1-1}_{|S_1}|$ can be estimated by
$$
|Df^{t_1-1}_{|S_1}|\ge C\frac{|f^{t_1}(K'')|}{|f(K'')|}
                        \ge C \frac{u_{n-1}}{u_n^\alpha}.
$$
\hfill\hfill\qed $\,\,$ (Lemma 3.5) 
\enddemo

\demo{Proof of Proposition 3.1} Assume first that $s=1$. 
This is an application of the previous Lemma 3.5. If 
$|R_n(U_{n+1})|\ge \frac{1}{10} u_n$ then
$$
|Df^{t_1}(f(x))|\ge C\cdot \frac{u_n}{u_{n+1}^\alpha}\cdot u_{n+1}^{\alpha-1}=
   C\cdot \frac{u_n}{u_{n+1}}\ge C\cdot \rho_n,
$$
where we used that $R_n(x)\notin U_{n+1}$.

\flushpar
In the other case when $|R_n(U_{n+1})|< \frac{1}{10} u_n$, we have 
$$
|Df^{t_1}(f(x))|\ge C\cdot \frac{u_{n-1}}{u_n^\alpha}\cdot u_n^{\alpha-1}
             =   C\cdot \frac{u_{n-1}}{u_{n}}\ge C\cdot \rho_n,
$$
where we used that in this case $f^{t_1-1}(x)\in R_n(U_{n+1})$
which is close to the boundary of $U_{n}$. The case with $s=1$ is finished. 

\bigskip

\flushpar
Now assume that $s\ge 2$. The proof will be split in two 
cases. Let $n\ge n_0\ge 1$ be big enough such that Lemma 3.3 and 3.4 can be 
applied. Let $K<\infty$ be the constant from Lemma 3.3 and $\theta<1$
the constant from Lemma 3.4.
 
\demo{Case I ( $f^T(0)=R_n^s(0)\in (-\frac34\cdot u_n,\frac34\cdot u_n)$)} Let 
$H\subset S_s$ be such that
$f^{T-1}=f^{t_s-1}:H\to (-\frac34 \cdot u_n, \frac34 \cdot u_n)$ is onto. 
This map has a monotone extension up to $U_n$. Hence it has a uniformly 
bounded distortion,
$$
\aligned
|Df^T(f(x))|&\ge C\cdot \frac{|f^{T-1}(H)|}{|H|}\cdot u_{n+1}^{\alpha-1}\\
            &\ge C\cdot \frac{u_n}{|S_s|}\cdot u_{n+1}^{\alpha-1}\\
            &\ge C\cdot \frac{u_n}
                             {\theta^{s-1}\cdot u_{n+1}^\alpha}
                  \cdot u_{n+1}^{\alpha-1}\\
            &\ge C\cdot \rho_n\cdot \theta^{-(s-1)}.
\endaligned
$$
\enddemo

\demo{Case II ( $f^T(0)=R_n^s(0)\notin (-\frac34\cdot u_n,\frac34\cdot u_n)$)} 
If the 
distortion of $f^{T-1}:S_s\to U_n$ is bounded by $K$ then we can give the same
argument as in case I:
$$
\aligned
|Df^T(f(x))|&\ge C\cdot \frac{|U_n|}{|S_s|}\cdot u_{n+1}^{\alpha-1}\\
            &\ge C\cdot \frac{u_n}
                             {\theta^{s-1}\cdot u_{n+1}^\alpha}
                  \cdot u_{n+1}^{\alpha-1}\\
            &\ge C\cdot \rho_n\cdot \theta^{-(s-1)}.
\endaligned
$$
Now let us assume that this distortion is bigger than $K$. Apply 
Lemma 3.3, which states $I_s\subset (-\frac34\cdot u_n,\frac34\cdot u_n)$. Then
$$
\aligned
|Df^{T}(f(0))|&=|Df^{t_{s-1}}(f(0))|\cdot |Df^{T-t_{s-1}}(f^{t_{s-1}}(f(0))|\\
              &\ge C\cdot \rho_n\cdot \theta^{-(s-2)}\cdot
                   |Df^{T-t_{s-1}}(f^{t_{s-1}}(f(0))|.
\endaligned
$$
For $s-1\ge 1$ we get this estimate from case I: 
$R_n^{s-1}(0)\in I_s\subset (-\frac34 \cdot u_n,\frac34\cdot u_n)$. 
When $s-1=1$ it follows from the Proof of Proposition 3.1 for $s=1$. 

\flushpar
The last factor can be
estimated by using the fact that $f^{T-t_{s-1}-1}:f(I_s)\to U_n$ has a monotone
extension up to $U_{n-1}$, see Lemma 2.2. It is essentially linear and its
derivative can be estimated
$$
|Df^{T-t_{s-1}-1}|_{f(I_s)}|\ge C\cdot \frac{|U_n|}{|f(I_s)|}
                            \ge C\frac{|U_n|}{\epsilon\cdot |f(U_n)|},
$$
where $\epsilon>0$ is given by Lemma 2.3. By taking $n_0\ge 1$ big enough we
can assume that $\epsilon>0$ is arbitrarily small.

\flushpar
Observe
that $|Df(f^{T-1}(f(0))|\ge C\cdot u_{n}^{\alpha-1}$ which follows from the 
fact that $|f^{T-1}(f(0))|\ge \frac34\cdot u_n$. 
We can finish the estimate for $|Df^{T}(f(0))|$ by observing that
$$
\aligned
|Df^{T-t_{s-1}}(f^{t_{s-1}}(f(0))|
&= 
|Df^{T-t_{s-1}-1}(f^{t_{s-1}}(f(0))|\cdot |Df(f^{T-1}(f(0))|\\
&\ge C\cdot \frac{|U_n|}{\epsilon\cdot |f(U_n)|}\cdot u_n^{\alpha-1}\\
&\ge C\cdot \frac{u_n}{\epsilon \cdot u_n^\alpha}\cdot u_n^{\alpha-1}
\ge C\cdot \frac{1}{\epsilon}\ge \frac{1}{\theta}.
\endaligned
$$	 
\hfill\hfill\qed $\,\,$ (Proposition 3.1)
\enddemo
\enddemo
\bigskip
\centerline{\bf 4. Telemann Decomposition of the Critical Orbit}
\bigskip

\flushpar
In this section we will prove Theorem B. Let $f$ be a unimodal
map such that 
$$
\sum_{n\ge 1} \sigma_n^{1/\alpha}<\infty.
$$
The existence of an absolutely continuous invariant probability
measure follows from [NS] in where it was shown that  Summability 
Condition on Derivatives
$$
\sum_{k\ge 1}|Df^k(f(0))|^{-\frac1\alpha}<\infty
$$
is 
sufficient for the existence of absolutely continuous invariant 
probability measures.

\flushpar
In the sequel we will prove that the summability of scaling factors
implies the Summability Condition of derivatives. 
Choose $n_0\ge 1$ big enough such that Proposition 3.1 can be applied
and moreover
$$
a=\sum_{{n\ge n_0} \atop {s\ge 0}}\frac{\theta^{\frac1\alpha\cdot s}}
                                     {(C\cdot\rho_n)^{\frac{1}{\alpha}}}= 
\frac{1}{1-\theta^{\frac1\alpha}} \sum_{n\ge n_0} 
\frac{1}{(C\cdot \rho_n)^\frac1\alpha}< 1,
$$
where $C$ and $\theta$ are the constants from Proposition 3.1.

\flushpar
Fix $k\ge 1$. The estimate for $|Df^k(f(0))|$ is based on the {\it Telemann 
decomposition} of the critical orbit up to time $k$. Consider the 
orbit of $f(0)$ up to time $k-1$. Let $m\ge 0$
be such that $U_{n_0+m}$ is the smallest central interval
which is crossed by this piece of the orbit:
$$
n_0+m=\max\{j\ge 0| \exists 0<i\le k \text{  }f^i(0)\in U_j\}
$$
and the last moment of crossing is denoted by
$$
k_m=\max\{1\le j\le k| f^j(0)\in U_{n_0+m}\}.
$$
The moments $k_m\le k_{m-1}\le \dots k_1\le k_0$ are such that $k_i$
is the last moment that the orbit hits $U_{n_0+i}$:
if $\{k_i<j\le k|f^j(0)\in U_{n_0+i-1}\}=\emptyset$ then $k_{i-1}=k_i$
otherwise
$$
k_{i-1}=\max\{k_i<j\le k|f^j(0)\in U_{n_0+i-1}\} \text{ with } 1\le i\le m.
$$
Let $r(k)=k-k_0$ and if $k_{i-1}\ne k_i$ then
$$
s_{i-1}(k)=\#\{k_i<j\le k_{i-1}|f^j(0)\in U_{n_0+i-1}\} \text{,  } 
1\le i\le m,
$$
the number of returns trough $U_{n_0+i-1}$.

\bigskip

\flushpar
The chain-rule applied to $Df^k(f(0))$ gives
$$
Df^k(f(0))=Df^{r(k)}(f^{k_0}(f(0)))\cdot
Df^{k_m}(f(0))\cdot
\prod_{i=0}^{m-1} Df^{k_i-k_{i+1}}(f^{k_{i+1}}(f(0))).
$$
The first factor can be estimated by using 

\proclaim{Proposition 4.1 ([G],[Ma])} Given $n_0\ge 1$ 
there exist constants $C>0$ and 
$\lambda>1$ such that 
$$
|Df^r(x)|\ge C\lambda^r,
$$
whenever $f^i(x)\notin U_{n_0}$ with $i\le r$.
\endproclaim

\flushpar
The other factors can be estimated by Proposition 3.1. 
The decomposition was 
set up to make Proposition 3.1 applicable to the factors:
$$
\aligned
|Df^k(f(0))|^{-\frac1\alpha}&
\le  \{ C\lambda^{r(k)}\cdot 
\prod_{ {i\le m} \atop  {k_i\ne k_{i+1}} }^{m} 
C\cdot \rho_{i}\cdot \theta^{-(s_i(k)-1)} \}^{-\frac1\alpha}\\
&\le  C(\frac{1}{\lambda^\frac1\alpha})^{r(k)}\cdot 
\prod_{ {i\le m} \atop  {k_i\ne k_{i+1}} }^{m} 
\frac{ (\theta^{\frac1\alpha})^{s_i(k)-1}}
     {(C\cdot \rho_{i})^{\frac1\alpha}}\\
\endaligned
$$

\proclaim{Lemma 4.2} Let $s_i,r$ and $s'_i,r'$ correspond to the Teleman 
decomposition of respectively $k$ and $k'$. If $k\ne k'$ then $r\ne r'$ or
$s_i\ne s'_i$ for some $i\ge 0$.
\endproclaim

\demo{Proof} Assume that $r=r'$ and $s_i=s'_i$ for all $i\ge 1$. We have to 
show 
that $k=k'$. Observe that $f^{k_m}(0)=R_{n_0+m}^{s_m}(0)$ but also
$f^{k'_m}(0)=R_{n_0+m}^{s'_m}(0)=R_{n_0+m}^{s_m}(0)$. So $k_m=k'_m$. Now
repeat this argument to show that $k_i=k'_i$ for $0\le i\le m$. In particular
we get $k_0=k'_0$. So
$$
k'=k'_0+r'=k_0+r=k.
$$
\hfill\hfill\qed $\,\,$ (Lemma 4.2)
\enddemo

\bigskip
\centerline{\bf Proof of the Summability Condition for Derivatives}
\bigskip

\flushpar
The number $a<1$ was defined in the beginning of this section. Let 
$\beta=\frac1\alpha$. 
$$
\aligned
\sum_{k\ge 0} |Df^k(f(0))|^{-\frac1\alpha}
&=\sum_{r\ge 0}\sum_{{k\ge 0}\atop {r(k)=r}  } |Df^k(f(0))|^{-\frac1\alpha}\\
&\le \sum_{r\ge 0}  C(\frac{1}{\lambda^\beta})^{r}\cdot
  \sum_{{k\ge 0} \atop {r(k)=r}} 
\prod_{{0\le i\le m} \atop {k_i\ne k_{i+1}} } 
\frac{ (\theta^{\beta})^{s_i(k)-1}}{(C\cdot \rho_{i})^{\beta}}.
\endaligned
$$
Now observe that for each $r$ there are no two  products appearing in the 
second sum which are formed with the same factors, Lemma 4.2. 
The sum of all possible different products
can be estimated by $1+a+a^2+a^3+\dots$. Hence  
$$
\aligned
\sum_{k\ge 0} |Df^k(f(0))|^{-\frac1\alpha}
&\le \sum_{r\ge 0}  C(\frac{1}{\lambda^\beta})^{r}\cdot 
(1+a+a^2+a^3+\dots)\\
&\le \frac{1}{1-a} \sum_{r\ge 0}  C(\frac{1}{\lambda^\beta})^{r}\\
&\le C\cdot \frac{1}{1-a} \cdot \frac{1}{1-\frac{1}{\lambda^{\beta}}}
<\infty.
\endaligned
$$

\bigskip
\centerline{\bf 5. Appendix}
\bigskip

\flushpar
In this appendix some basic notions of interval dynamics are collected. The
details can be found in [MS].

\bigskip

\flushpar
The hyperbolic length of an interval $I\subset T\subset [-1,1]$ within $T$ 
is defined to be
$$
\text{hyp}(I,T)=\ln\frac{|L\cup I|\cdot |R\cup I|}{|L|\cdot |R|},
$$
where $L,R\subset T$ are the connected components of $T\setminus I$ and 
$|J|$ stands for the length of the interval $J\subset [-1,1]$.

\flushpar
The Schwarzian derivative of a $C^3$ map $f:[-1,1]\to [-1,1]$ is
$$
Sf(x)=\frac{D^3f(x)}{Df(x)}-\frac32\cdot\frac{D^2f(x)}{Df(x)}.
$$
where $D^if(x)$ stands for the $i^{th}$ derivative of $f$ in $x\in [-1,1]$.

\proclaim{Expansion-Lemma} If $f:[-1,1]\to [-1,1]$ has $Sf(x)\le 0$ for all
$x\in [-1,1]$ and $f^n|T$ is monotone then
$$
\text{hyp}(f^n(I),f^n(T))\ge \text{hyp}(I,T),
$$
where $I\subset T$.
\endproclaim

\proclaim{Koebe-Distortion-Lemma} For each $\tau>0$ there exists 
$1\le K(\tau)<\infty$ with the following property. 
Let $f^n:T\to f^n(T)$ be monotone and $Sf(x)\le 0$
for all $x\in [-1,1]$. If $I\subset T$ is an interval such that both component
$L,R\subset T\setminus I$ satisfy
$$
\frac{|f^n(L)|}{|f^n(T)|},\frac{|f^n(R)|}{|f^n(T)|}\ge \tau
$$
then $f^n|I$ has bounded distortion
$$
\frac{|Df^n(x)|}{|Df^n(y)|}\le K(\tau),
$$
for all $x,y\in I$.
Moreover $K(\tau)\to 1$ when $\tau\to \infty$.
\endproclaim

\flushpar
An $S-$unimodal map is a $C^3$ endomorphism $f:[-1,1]\to [-1,1]$ such that

\parindent=15pt
\item{1)} $f(\pm1)=-1$,
\item{2)} $Df(x)>0$ for $x<0$,
\item{3)} $Df(x)<0$ for $x>0$,
\item{4)} $Df(0)=0$ and up to a $C^1$ change of coordinates
$f$ equals locally $x\to |x|^{\alpha}$ with $\alpha>1$. 
The point $x=0$ is called the critical point and $\alpha>1$ is called the
critical exponent of $f$. 
\item{5)} $Sf(x)<0$, $x\ne 0$.

\flushpar
The orbit of the critical point $x=0$ is called the critical orbit. The 
critical orbit is said to be recurrent if it accumulates at the critical 
point.

\bigskip

\flushpar
\demo{Usage of constants} Every uniform constant, that is a constant which
is independent of the actual map, appearing in estimates will be denoted by 
$C$. Constants which play a specific role in the statement of Lemmas will
usually have a specific name.
\enddemo
\bigskip
\centerline{\bf References}
\bigskip

\parindent=25pt
\item{[BKNS]} H.Bruin, G.Keller, T.Nowicki, S.van Strien, {\it Wild
            Attractors Exist}, to appear in Ann. of Math.

\item{[B]} L.A.Bunimovich, {\it A Transformation of the Circle},
           Math.Notes {\bf 8} no. 2 (1970), 587-592.

\item{[Br]} H.Bruin, {\it Invariant Measures of Interval Maps},
            Thesis (1994), Technical University Delft, Delft.

\item{[BL1]} A.M.Blokh, M.Lyubich, {\it Attractors of Maps on the interval}
            Func.Anal. and Appl. {\bf 21} (1987), 32-46.

\item{[BL2]} A.M.Blokh, M.Lyubich, {\it Measure of Solenoidal Attractors
            of Unimodal Maps of the Segment} Math.Notes {\bf 48} no. 5 (1990),
            1085-1090.

\item{[G]} J.Guckenheimer, {\it Limit Sets of $S-$unimodal Maps with
            Entropy Zero}, Commun.Math.Phys. {\bf 110} (1987), 655-659.

\item{[GJ]} J.Guckenheimer, S.Johnson, {\it Distortion of $S-$unimodal Maps},
            Ann. Math. {\bf 132} (1990), 71-130. 

\item{[GS]} J.Graczyk, G.Swiatek, {\it Hyperbolicity in the Real
            Quadratic Family}, Report no. PM 192 PennState (1995).

\item{[HK]} F.Hofbauer, G.Keller, {\it Quadratic Maps without Asymptotic
            Measure} Commun.Math. Phys. {\bf 127} (1990), 319-337.

\item{[Ja]} M.V.Jakobson, {\it Absolutely Continuous Invariant Measures 
           for one-parameter Families of One-dimensional Maps} 
           Commun.Math.Phys. {\bf 81} (1981), 39-88.

\item{[Jo]} S.Johnson, {\it Singular Measures without restrictive interval},
           Commun.Math.Phys. {\bf 110} (1987), 185-190.

\item{[K]} G.Keller, {\it Exponents, Attractors and Hopf Decomposition
            for Interval Maps} Erg.Th. \& Dyn.Sys. {\bf 10} (1990), 717-744. 

\item{[L1]} M.Lyubich, {\it Combinatorics, Geometry and Attractors of
            Quasi-Quadratic Maps}, Ann. of Math {\bf 140} (1994), 347-404.

\item{[L2]} M.Lyubich, {\it Geometry of Quadratic Polynomials: Moduli,
            Rigidity, and Local Connectivity},
            IMS at Stony Brook preprint 1993/9. 

\item{[L3]} M.Lyubich, {\it Dynamics of Quadratic Polynomials, II 
            Rigidity}, IMS at Stony Brook preprint 1995/14. 

\item{[L4]} M.Lyubich, {\it Dynamics of Quadratic Polynomials, III 
            Parapuzzles}, IMS at Stony Brook preprint 1996/5.

\item{[LM]} M.Lyubich, J.W.Milnor, {\it The Fibonacci Unimodal Map},
            J.A.M.S. {\bf 6} (1993), 425-457.

\item{[Ld]} F.Ledrappier, {\it Some Properties of Absolutely Continuous
            Invariant Measures of an Interval}, Erg.Th. \& Dyn.Sys. 
            {\bf 1} (1981), 77-93. 

\item{[Ma]} R.Ma\~n\'e, {\it Hyperbolicity, Sinks, and Measure in 
            One-dimensional Dynamics}, Commun.Math.Phys. {\bf 100} (1985), 
            495-524 and Erratum. Commun.Math.Phys. {\bf 112} (1987), 721-724. 

\item{[M]} M.Martens, {\it Interval Dynamics}, Thesis (1990), 
            Technical University Delft, Delft. or
 
\item{[M$'$]} M.Martens, {\it Distortion Results and Invariant Cantor Sets of
           Unimodal Maps}, Erg.Th. \& Dyn.Sys. {\bf 14} (1994), 331-349.

\item{[Mi]} M.Misiurewicz, {\it Absolutely Continuous Measures for certain
            maps of an Interval}  Publ.Math. I.H.E.S. {\bf 53} (1981),
            17-51.
 
\item{[MS]} W.de Melo and S.van Strien, {\it One-dimensional Dynamics},
           Springer-Verlag, 1993.

\item{[NS]} T.Nowicki, S.van Strien, {\it Absolutely Continuous Measures under
           a Summability Condition}, Invent.Math. {\bf 93} (1988), 619-635.

\item{[NU]} J.von Neumann, S.M.Ulam, {\it On Combinatorics of Stochastic and
            Deterministic Processes}, Bull.A.M.S. {\bf 53} (1947), 1120.

\item{[R]}  D.Ruelle, {\it Applications Conservants une mesure absolutement 
            continu par rapport \`a $dx$ sur $[0,1]$}, Commun.Math.Phys. 
            {\bf 55} (1977), 47-51.

\item{[VT]} J.J.P.Veerman, F.M.Tangerman, {\it Scalings in Circle Maps (I)},
            Commun.Math.Phys. {\bf 134} (1990), 89-107.


\bye